\documentclass[a4paper,11pt]{article}
\usepackage{amsmath, amssymb, amsthm}
\usepackage{graphicx} 

\numberwithin{equation}{section}

\begin{document}

\def\R{{\mathbb R}} \def\PP{{\mathcal P}} \def\C{{\mathbb C}}
\def\HH{{\mathcal H}} \def\Markoc{{\rm Markov}} \def\Skein{{\rm Skein}}
\def\K{{\mathbb K}} \def\LL{{\mathcal L}} \def\BB{{\mathcal B}}
\def\Sym{{\rm Sym}} \def\Markov{{\rm Markov}} \def\Z{{\mathbb Z}}
\def\N{{\mathbb N}} \def\S{{\mathbb S}}

\title{\bf{HOMFLYPT Skein module of singular links}}
 
\author{\textsc{Luis Paris}}

\date{\today}

\maketitle

\begin{abstract}
\noindent
This paper is a presentation, where we compute the HOMFLYPT Skein module of singular links in the $3$-sphere. This calculation 
is based on some results previously proved by Rabenda and the author on Markov traces on singular Hecke algebras, 
as well as on classical techniques that allow to pass from the framework of Markov traces on Hecke algebras to the framework of HOMFLYPT Skein 
modules. Some open problems on singular Hecke algebras are also presented.
\end{abstract}




\section{Introduction}

Knot theory had a notable renewable in the 80's with the emergence of new knot invariants such as the Jones polynomial 
\cite{Jones1}, \cite{Jones2} and the HOMFLYPT polynomial \cite{HOMFLY}, \cite{PrzTra1}. The latest one is defined by the 
following theorem.

\bigskip\noindent
{\bf Theorem 1.1} (Freyd, Yetter, Hoste, Lickorish, Millett, Ocneanu \cite{HOMFLY}, Przytycki, Traczyk \cite{PrzTra1}).
{\it Let $\LL$ be the set of (isotopy classes) of oriented links in the sphere $\S^3$. Then there exists a unique invariant 
$I:\LL \to \C(t,x)$ which is $1$ on the trivial knot, and which satisfies the relation
\begin{equation}
t^{-1} \cdot I(L_+) -t \cdot I(L_-) = x \cdot I(L_0)\,,
\label{eq11}
\end{equation}
for all links $L_+,L_-,L_0 \in \LL$ that have the same link diagram except in the neighborhood of a crossing where they are 
like in Figure 1.1.}
\begin{figure}[htb]
\bigskip
\centerline{
\setlength{\unitlength}{0.5cm}
\begin{picture}(9,2.5)
\put(0,0.5){\includegraphics[width=4.5cm]{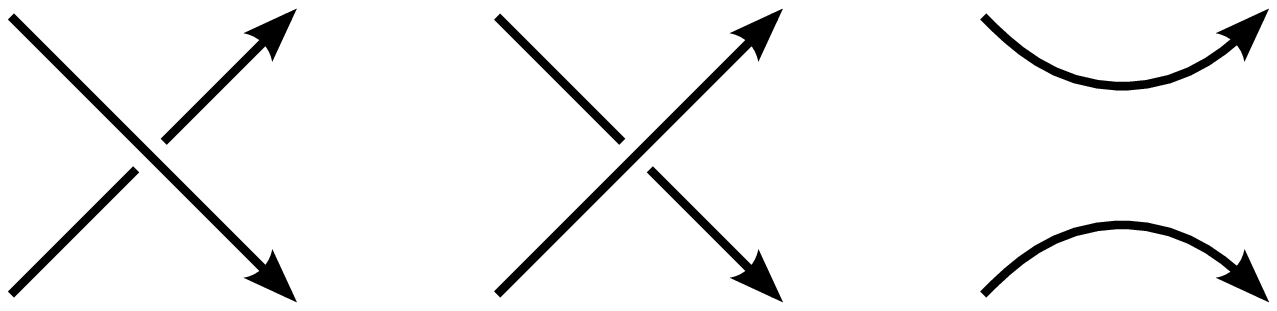}}
\put(0.8,-0.1){\small  $L_+$}
\put(4.3,-0.1){\small  $L_-$}
\put(7.8,-0.1){\small  $L_0$}
\end{picture}}
\bigskip
\centerline{{\bf Figure 1.1.} The links $L_+$, $L_-$, and $L_0$.}
\end{figure}

\bigskip\noindent
Since then, knot theorists wonder about possible extensions of this result to other sets of like-knots such as the set of links in a
3-manifold, or the set of singular links in the $3$-sphere.

\bigskip\noindent
Recall that a {\it singular link} on $n$ components is defined to be an immersion of $n$ circles in the sphere $\S^3$ which admits 
only finitely many singularities that are all ordinary double points. By \cite{Kauff1}, two singular link diagrams 
represent the same singular link (up to isotopy)
if and only if one can pass from one to the other by a finite sequence of  ordinary or singular Reidemeister moves (see 
Figures 1.2 and 1.3).
\begin{figure}[htb]
\bigskip
\centerline{
\setlength{\unitlength}{0.5cm}
\begin{picture}(26,12)
\put(0,0){\includegraphics[width=13cm]{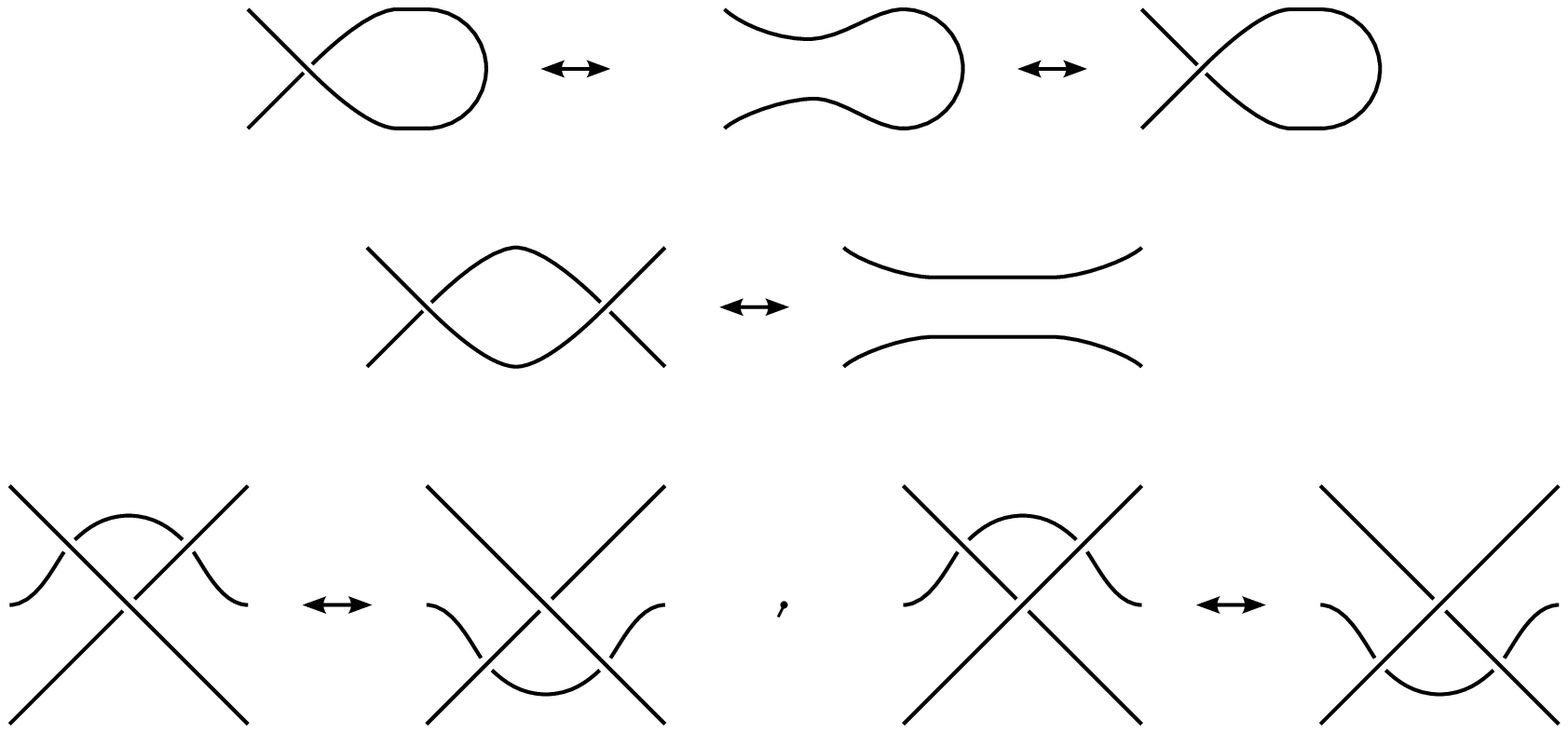}}
\end{picture}}
\bigskip
\centerline{{\bf Figure 1.2.} Ordinary Reidemeister moves.}
\end{figure}
\begin{figure}[htb]
\bigskip
\centerline{
\setlength{\unitlength}{0.5cm}
\begin{picture}(30,8)
\put(0,0){\includegraphics[width=15cm]{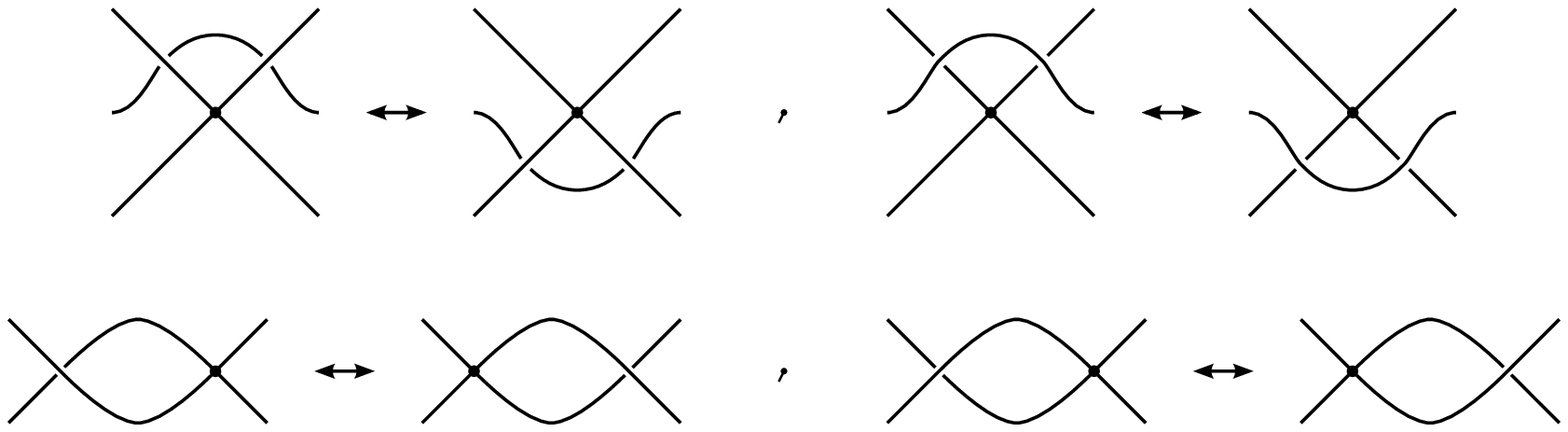}}
\end{picture}}
\bigskip
\centerline{{\bf Figure 1.3.} Singular Reidemeister moves.}
\end{figure}

\bigskip\noindent
Let $\LL$ be a set of like-knots. We say 
that an invariant $I:\LL \to \C(t,x)$ satisfies the {\it HOMFLYPT Skein relation} if the relation \eqref{eq11} holds for all 
links $L_+,L_-,L_0 \in \LL$ that have the same link diagram except in the neighborhood of a crossing where they are like in Figure 
1.1. It has been quickly observed that, in general, there are many invariants that satisfy the HOMFLYPT Skein relation and  
that are $1$ on the trivial knot. However, the condition that the invariant is $1$ on the trivial knot is secondary, and, moreover, one 
can view the set of invariants that satisfy the HOMFLYPT Skein relation as a vector space over $\C(x,t)$.
So, the general question is in fact to determine this vector space.

\bigskip\noindent
Let $\LL$ be a set of like-knots. Define the {\it HOMFLYPT Skein module} of $\LL$, denoted by $\Skein(\LL)$, to be the quotient of 
the vector space $\C(x,t)[\LL]$ freely spanned by $\LL$, by the relations
\[
t^{-1} \cdot L_+ - t \cdot L_- = x \cdot L_0\,,
\]
for all links $L_+,L_-,L_0 \in \LL$ that have the same link diagram except in the neighborhood of a crossing where they are like in 
Figure 1.1. Note that the space of invariants of $\LL$ that satisfy the HOMFLYPT Skein relation is the space of linear forms 
on $\Skein(\LL)$.

\bigskip\noindent
The Skein module was calculated for the set of links in a solid torus by Hoste, Kidwell \cite{HosKid1} and, independently, Turaev \cite{Turae1}. They result was extended by Przytycki \cite{Przyt1} to the set of links in the direct product $F \times I$ of a surface $F$ with the interval. In this case, $\Skein(\LL)$ can be endowed with a structure of algebra. The product of two links $L_1$ and $L_2$ (modulo the Skein relations) is the link obtained placing $L_2$ above $L_1$. Note that the Skein module of singular links can be also endowed with a structure of algebra following the same rules.

\bigskip\noindent
{\bf Theorem 1.2} (Przytycki \cite{Przyt1}).
{\it Let $\LL$ be the set of links in the direct product $F \times I$ of a surface $F$ with the interval $I$. Then 
$\Skein(\LL)$ is isomorphic to the symmetric algebra $S\C(t,x)[\hat \pi^0]$ on the vector space $\C(t,x)[\hat \pi^0]$ freely 
spanned by the set $\hat\pi^0$ of conjugacy classes of nontrivial elements of $\pi_1(F)$.}

\bigskip\noindent
The purpose of this paper is to present an approach to the calculation of HOMFLYPT Skein modules via the study of different 
sorts of braid groups and monoids and their associated generalized Hecke algebras. This will be done through the study of a 
particular example: the singular links in the $3$-sphere. However, the ideas presented here can be easily extended to other cases. 
In particular, a careful reading of \cite{Lambr1} shows how to use these techniques to calculate the HOMFLYPT Skein module of 
the solid torus.

\bigskip\noindent
The main result of this paper is:

\bigskip\noindent
{\bf Theorem 1.3.}
{\it Let $\LL$ be the set of oriented singular links in the sphere $\S^3$. Then $\Skein(\LL)$ is isomorphic to the polynomial 
algebra $\C(x,t)[\hat X,\hat Y]$ in the two variables $\hat X$ and $\hat Y$, where $\hat X$ and $\hat Y$ are represented by 
the links $L_X$ and $L_Y$ drawn in Figure 1.4}
\begin{figure}[htb]
\bigskip
\centerline{
\setlength{\unitlength}{0.5cm}
\begin{picture}(11,3.5)
\put(0,0.5){\includegraphics[width=5.5cm]{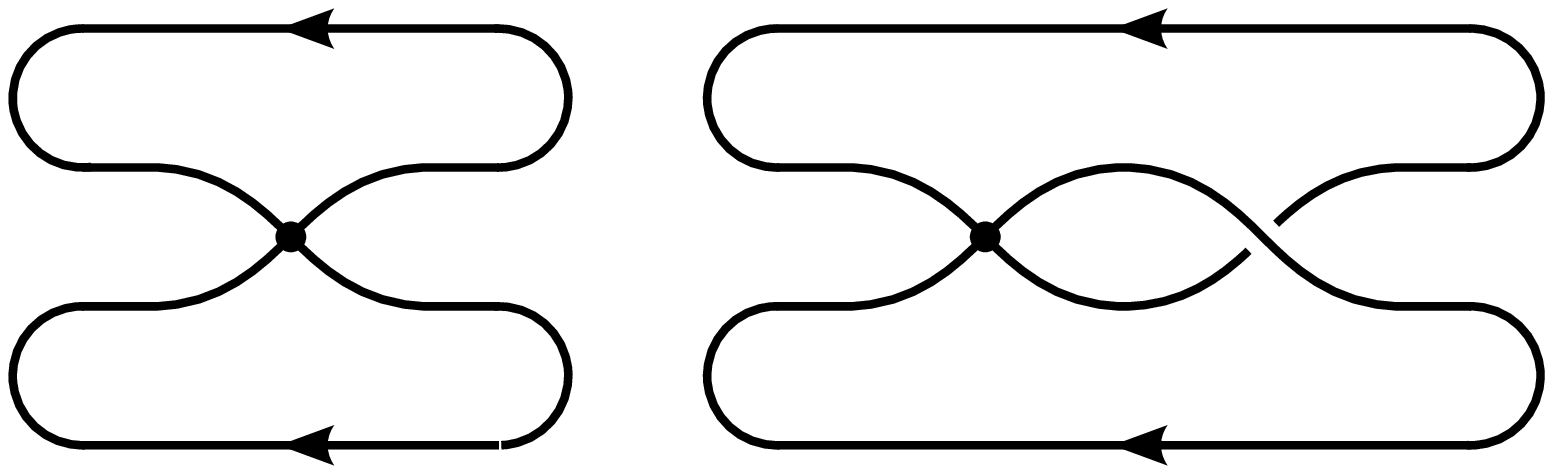}}
\put(1.8,-0.3){\small $L_X$}
\put(7.8,-0.3){\small $L_Y$}
\end{picture}}
\bigskip
\centerline{{\bf Figure 1.4.} Generators of the Skein module of singular links.}
\end{figure}

\bigskip\noindent
The proof of Theorem 1.3 consists essentially in translating the main result of \cite{ParRab1}, which concerns Markov traces on singular Hecke 
algebras, in terms of HOMFLYPT Skein modules. On the other hand, some open questions will be presented along the text, and the 
proof of Theorem~1.3 will also serve as a pretext to present them.


\section{Markov module and HOMFLYPT Skein module}

Let $\PP=\{P_1,\dots, P_n\}$  be a set of  $n$ distinct punctures in the plane $\R^2$ (except mention of the contrary, we will always 
assume $P_k=(k,0)$ for all $1\le k \le n$). A {\it singular braid} on $n$ strands is defined to be an $n$-tuple $\beta=(b_1, 
\dots, b_n)$ of smooth paths, $b_k:[0,1] \to\R^2 \times [0,1]$, such that
\begin{itemize}
\item
there exists a permutation  $\chi \in \Sym_n$ such that $b_k(0)=(P_k,0)$ and $b_k(1)=(P_{\chi (k)}, 1)$ for all $1 \le k\le n$;
\item
$b_k(t)$ runs monotonically on the second coordinate for all  $1 \le k\le n$;
\item
the image of $b_1 \sqcup \cdots \sqcup b_n$ has finitely many singularities (called {\it singular points}), that are all ordinary 
double points.
\end{itemize}

\bigskip\noindent
The isotopy classes of singular braids form a monoid called {\it singular braid monoid} (on $n$ strands) and denoted by $SB_n$. 
The multiplication in this monoid is the concatenation of (singular) braids.

\bigskip\noindent
{\bf Theorem 2.1} (Baez \cite{Baez1}, Birman \cite{Birma1}). 
{\it The monoid $SB_n$ has a monoid presentation with generators
\[
\sigma_1, \dots, \sigma_{n-1}, \sigma_1^{-1}, \dots, \sigma_{n-1}^{-1}, \tau_1, \dots, \tau_{n-1}\,,
\]
and relations 
\[\begin{array}{cl}
\sigma_k \sigma_k^{-1} = \sigma_k^{-1} \sigma_k = 1 &\quad\text{for } 1 \le k \le n-1\,,\\
\sigma_k \tau_k = \tau_k \sigma_k &\quad\text{for } 1 \le k \le n-1\,,\\
\sigma_k \sigma_l \sigma_k = \sigma_l \sigma_k \sigma_l\,, &\quad\text{if } |k-l|=1\,,\\
\sigma_k \sigma_l \tau_k = \tau_l \sigma_k \sigma_l &\quad\text{if } |k-l|=1\,,\\
\sigma_k \sigma_l = \sigma_l \sigma_k &\quad \text{if } |k-l| \ge 2\,,\\
\sigma_k \tau_l = \tau_l \sigma_k &\quad \text{if } |k-l| \ge 2\,,\\
\tau_k \tau_l = \tau_l \tau_k &\quad \text{if } |k-l| \ge 2\,. 
\end{array}\]}

\bigskip\noindent
The braid $\sigma_k$ in the above theorem is the standard $k$-th generator of the braid group $B_n$ (see Figure 2.1). The braid 
$\tau_k$ is a singular braid with a unique singular crossing between the $k$-th strand and the $(k+1)$-th strand (see Figure 
2.1). 
\begin{figure}[htb]
\bigskip
\centerline{
\setlength{\unitlength}{0.5cm}
\begin{picture}(13,7)
\put(1,0){\includegraphics[width=6cm]{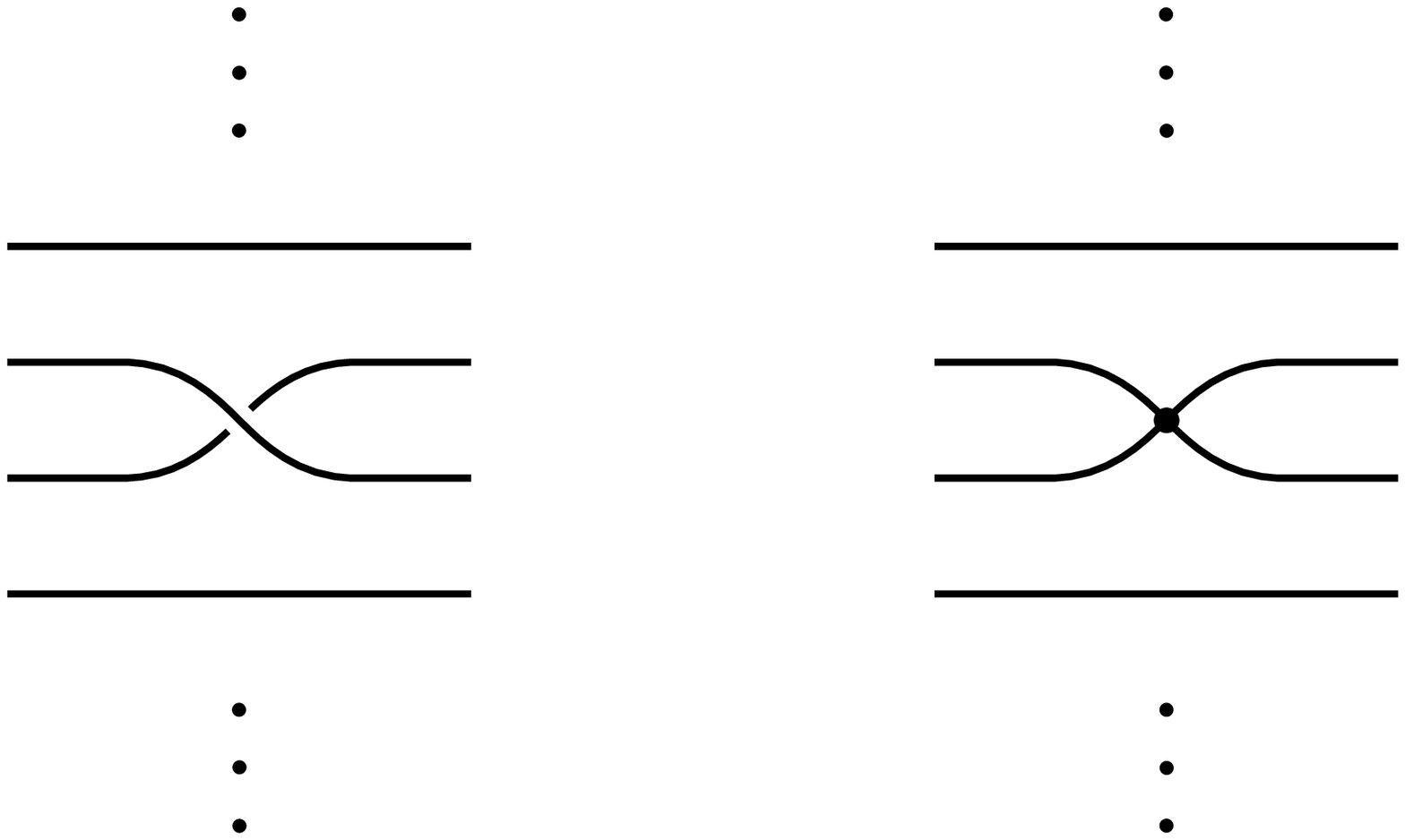}}
\put(0.5,2.9){\small $k$}
\put(-0.6,3.9){\small $k+1$}
\put(-1.9,3.4){\small $\sigma_k =$ }
\put(8.5,2.9){\small $k$}
\put(7.4,3.9){\small $k+1$}
\put(6.2,3.4){\small $\tau_k =$ }
\end{picture}}
\bigskip
\centerline{{\bf Figure 2.1.} Generators of $SB_n$.}
\end{figure}

\bigskip\noindent
From a singular braid $\beta$ we can construct a singular link connecting 
the point  $(P_k,1)$ to the point $(P_k,0)$ for all $1 \le k \le n$
(see Figure 2.2). This link is denoted by $\hat \beta$ and is called the {\it closure} of $\beta$. By \cite{Birma1},
every singular link is a closed singular braid. 
\begin{figure}[htb]
\bigskip
\centerline{
\setlength{\unitlength}{0.5cm}
\begin{picture}(22,9)
\put(1,0){\includegraphics[width=10.5cm]{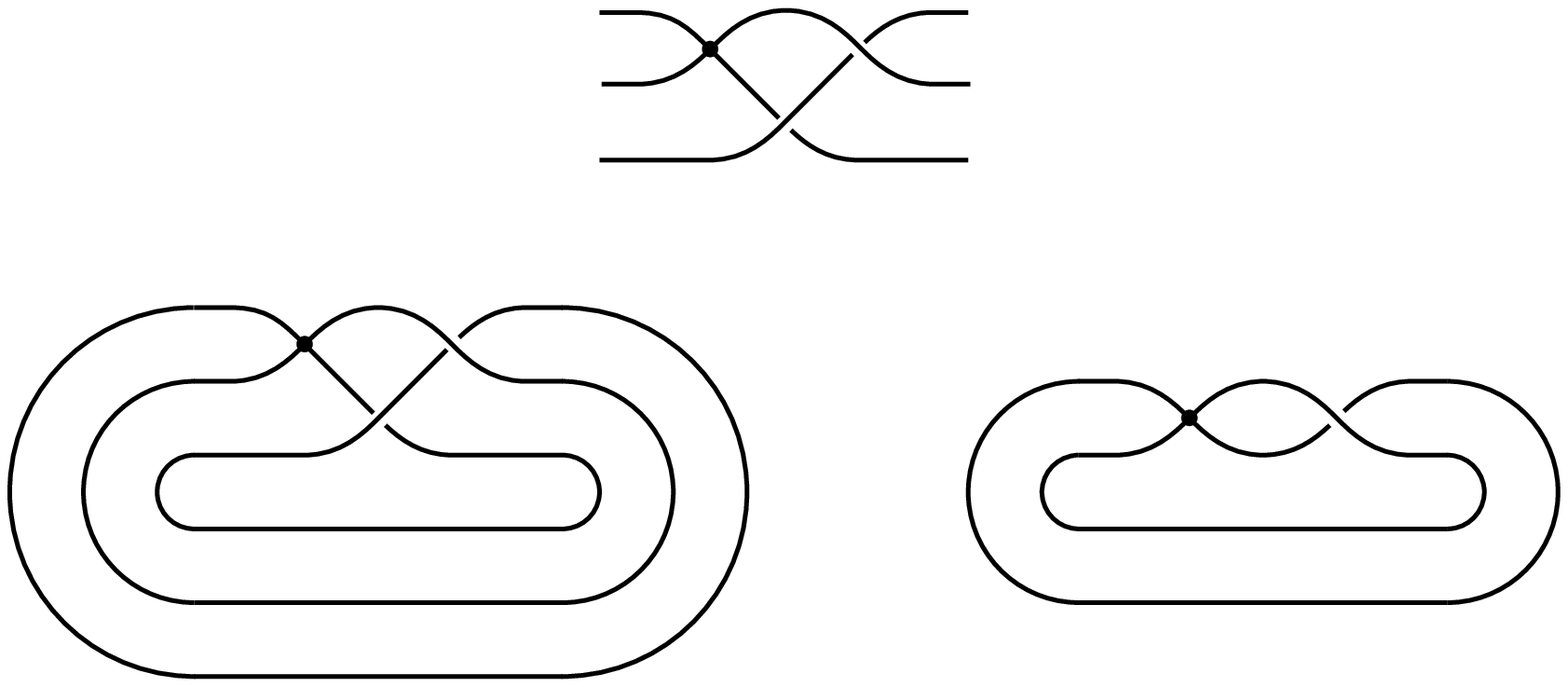}}
\put(-0.6,2.4){\small $\hat \beta=$}
\put(7.5,8){\small $\beta=$}
\put(12.4,2.4){$=$}
\end{picture}}
\bigskip
\centerline{{\bf Figure 2.2.} A closed braid.}
\end{figure}

\bigskip\noindent
We denote by $\sqcup SB= \sqcup_{n=1}^\infty SB_n$ the disjoint union of all singular braid monoids. We use the notation 
$(\beta,n)$ to denote a singular braid $\beta$ in $SB_n$  if we need to specify the number $n$ of strands.

\bigskip\noindent
Two singular braids $(\alpha,n)$ and $(\beta,m)$ are said to be connected by a {\it Markov move} if either
\begin{itemize}
\item
$n=m$ and there exist $\gamma_1, \gamma_2 \in SB_n$ such that $\alpha= \gamma_1 \gamma_2$ and $\beta = \gamma_2 \gamma_1$; or
\item
$m=n+1$ and $\beta = \alpha \sigma_n^{\pm 1}$; or
\item
$n=m+1$ and $\alpha = \beta \sigma_m^{\pm 1}$.
\end{itemize}

\bigskip\noindent
{\bf Theorem 2.2} (Gemein \cite{Gemei1}). 
{\it Let $(\alpha,n)$, $(\beta,m)$ be two singular braids. Then $\hat \alpha$ and $\hat \beta$ are isotopic if and only if $(\alpha,n)$ and $(\beta,m)$ are connected by a finite sequence of Markov moves. }

\bigskip\noindent
We turn now to apply this theorem to obtain a version of the HOMFLYPT Skein module of singular links in terms of singular Hecke algebras.

\bigskip\noindent
The {\it singular Hecke algebra}, denoted by  $\HH (SB_n)$, is defined to be the quotient of the monoid algebra $\C(q)[SB_n]$ 
by the relations
\begin{equation}
\sigma_k^2 = (q-1) \sigma_k +q\,, \quad 1 \le k \le n-1\,.
\label{eq21}
\end{equation}

\bigskip\noindent
Note that the singular Hecke algebra is an infinite dimensional $\C(q)$-vector space (except for $n=1$). However, it can be 
endowed with a graduation, and each term of the graduation is of finite dimension (see \cite{ParRab1}). This graduation is 
defined as follows.

\bigskip\noindent
For $n\ge 2$ and $d\ge 0$, we denote by $S_dB_n$ the set of singular braids with $n$ strands and $d$ singular points, and we 
denote by $\C(q)[S_dB_n]$ the subspace of $C(q)[SB_n]$ spanned by $S_dB_n$. Note that  $S_0B_n$ is the braid group $B_n$ on $n$ 
strands, and $\C(q)[S_0B_n] = \C(q)[B_n]$ is the group algebra of $B_n$. The monoid algebra $\C(q)[SB_n]$ is naturally graded by
\[
\C(q) [SB_n] = \bigoplus_{d=0}^{+ \infty} \C(q) [S_dB_n]\,.
\]
Now, the relations \eqref{eq21} that define the singular Hecke algebra are all homogeneous (of degree $0$), thus the graduation 
of $\C(q)[SB_n]$ induces a graduation on $\HH(SB_n)$,
\[
\HH (SB_n) = \bigoplus_{d=0}^{+ \infty} \HH (S_dB_n)\,,
\]
where $\HH(S_dB_n)$ is the subspace of $\HH(SB_n)$ spanned by $S_dB_n$.

\bigskip\noindent
Several elementary questions on singular Hecke algebras are still open. Here are two of them.

\bigskip\noindent
{\bf Question 2.3.}  Note that $\HH(S_0B_n) = \HH(B_n)$ is the Hecke algebra of the symmetric group, thus  
$\HH(S_dB_n)$ is a representation of $\HH(B_n)$. It would be interesting to characterize this representation. Actually, the 
dimension  itself (over $\C(q)$) of $\HH(S_dB_n)$ is unknown, even for $d=2$.

\bigskip\noindent
{\bf Question 2.4.} The natural inclusion $SB_n \hookrightarrow SB_{n+1}$ induces a homomorphism $\iota_n: \HH(SB_n) \to \HH(SB_{n+1})$. We do not know whether $\iota_n$ is injective.

\bigskip\noindent
Now, we introduce a new variable $z$, 
we set $\HH_z (SB_n) = \C (z,q) \otimes_{\C(q)} \HH(SB_n)$ for all $n \ge 1$,
and we consider the direct sum $\oplus_{n=1}^\infty \HH_z(SB_n)$. Like for 
the singular braids, we use the notation $(a,n)$ to denote an element $a \in \HH_z(SB_n)$ if we need to 
specify the number $n$ of strands. 

\bigskip\noindent
The {\it Markov module} of $\sqcup SB$, denoted by $\Markov(\sqcup SB)$, is defined to be the quotient of the space 
$\oplus_{n=1}^\infty \HH_z(SB_n)$ by the relations
\begin{itemize}
\item
$(ab,n) = (ba,n)$ for all $n \ge 1$ and all $a,b \in \HH_z(SB_n)$;
\item
$(a,n) = (\iota_n(a), n+1)$ for all $n \ge 1$ and all $a \in \HH_z(SB_n)$;
\item
$(\iota_n(a) \sigma_n,n+1) = z \cdot (a,n)$ for all $n \ge 1$ and all $a \in \HH_z(SB_n)$.
\end{itemize}

\bigskip\noindent
The space $\Markov(\sqcup SB)$ can be endowed with a structure of $\C(z,q)$-algebra as follows.  Let  
$[\alpha,n]$ denote the element of  $\Markov(\sqcup SB)$ represented by a braid $(\alpha,n)$. Let $(\alpha,n)$ and $(\beta,m)$ be two 
braids. Then the product $[\alpha,n] \cdot [\beta,m]$ is represented by the braid in $SB_{n+m}$ obtained placing  $\beta$ above 
$\alpha$. Note that the unit for this multiplication is represented by the trivial braid in $SB_1=B_1=\{1\}$.

\bigskip\noindent
{\bf Lemma 2.5.} {\it The above defined multiplication in  $\Markov(\sqcup SB)$ is commutative.}

\bigskip\noindent
{\bf Proof.}
Let $(\alpha,n)$, $(\beta,m)$ be two singular braids. Let $(\alpha \ast \beta,n+m)$ be the braid obtained placing $\beta$ above 
$\alpha$. So, $[\alpha,n] \cdot [\beta,m] = [\alpha \ast \beta,n+m]$. Let  $\sigma_{n,m} \in B_{n+m}$ be the braid pictured in 
Figure~2.3. Observe that $\sigma_{n,m}(\beta\ast\alpha) \sigma_{n,m}^{-1} = (\alpha \ast \beta)$, thus $[\alpha,n] \cdot 
[\beta,m] = [\beta,m] \cdot [\alpha,n]$.
\qed
\begin{figure}[htb]
\bigskip
\centerline{
\setlength{\unitlength}{0.5cm}
\begin{picture}(9,7)
\put(0,0){\includegraphics[width=4.5cm]{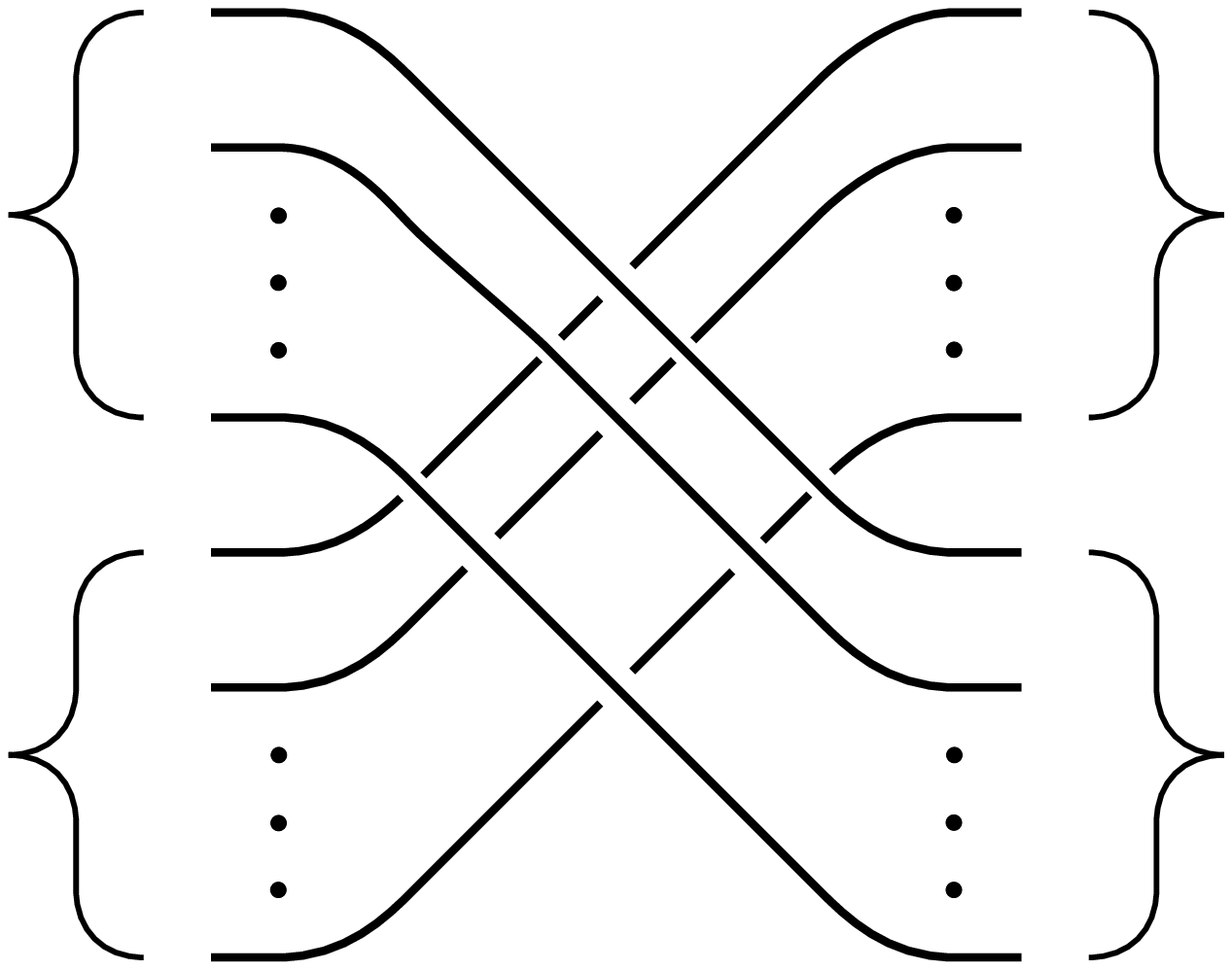}}
\put(-0.6,1.4){\small $n$}
\put(-0.7,5.4){\small $m$}
\put(9.2,1.4){\small $m$}
\put(9.2,5.4){\small $n$}
\end{picture}}
\bigskip
\centerline{{\bf Figure 2.3.} The braid $\sigma_{n,m}$.}
\end{figure}

\bigskip\noindent
Now, the link between the HOMFLYPT Skein module of singular links and the Markov module of singular braids is given by the 
following.

\bigskip\noindent
{\bf Theorem 2.6.} {\it 
Let $\LL$ be the set of singular links in the sphere $\S^3$. Set
\begin{gather*}
z = \frac{q-1}{1-qy} \quad \Leftrightarrow \quad y = \frac{z-q+1}{qz}\,,\\
t=\sqrt{yq} \,, \quad x=\sqrt{q} - \frac{1}{\sqrt{q}}\,.
\end{gather*}
Let $\K = \C(\sqrt{y},\sqrt{q})$. Then $\K \otimes \Markov (\sqcup SB)$ is isomorphic to $\K \otimes \Skein (\LL)$.}

\bigskip\noindent
{\bf Proof.}
In \cite{Jones1} Jones gives formulas to pass from Ocneanu's trace to the HOMFLYPT polynomial. In order to prove the above theorem, it suffices 
to slightly adapt these formulas to the context of the theorem.

\bigskip\noindent
For $(\beta,n) \in \sqcup SB$  we denote by $[\beta,n]$ the element of $\Markov(\sqcup SB)$ represented by  $(\beta,n)$. 
Similarly, for $L \in \LL$ we denote by $[L]$ the element of $\Skein(\LL)$ represented by $L$.

\bigskip\noindent
Let $ \psi_1 : \sqcup SB \to \K \otimes \Markov (\sqcup SB)$ be the map defined by
\[
\psi_1(\alpha,n) = \left( \frac{q-1}{1-qy} \right)^{-n+1} \left( \sqrt{y} \right)^{\varepsilon (\alpha) -n+1}
 [\alpha,n]\,,
\]
where $\varepsilon : SB_n \to \Z$ is the homomorphism defined by
\[
\varepsilon(\sigma_k) = 1\,,\ \varepsilon (\sigma_k^{-1}) = -1\,,\ \varepsilon(\tau_k)=0\,, \quad \text{for } 1 \le k\le
n-1\,.
\]
Let $(\alpha,n)$,  $(\beta,m)$ be two singular braids. We start showing that, if $\hat \alpha = \hat \beta$, then $\psi_1 
(\alpha,n) = \psi_1 (\beta,m)$. By Theorem 2.2, in order to do so, it suffices to consider the following three cases: 
\begin{enumerate}
\item
$n=m$ and there exist $\gamma_1, \gamma_2 \in SB_n$ such that $\alpha =\gamma_1 \gamma_2$ and $\beta = \gamma_2 \gamma_1$;
\item
$m=n+1$ and $\beta = \alpha \sigma_n$;
\item
$m=n+1$ and $\beta = \alpha \sigma_n^{-1}$.
\end{enumerate}

\bigskip\noindent
Suppose that $n=m$ and there exist $\gamma_1, \gamma_2 \in SB_n$ such that $\alpha= \gamma_1 \gamma_2$ and $\beta = \gamma_2 
\gamma_1$. By definition we have $[\alpha,n]=[\beta,n]$ and $\varepsilon(\alpha) = \varepsilon(\beta)$, thus $\psi_1(\alpha) = 
\psi_1(\beta)$. Suppose that $m=n+1$ and $\beta = \alpha \sigma_n$. Then 
\[\begin{array}{lcl}
\psi_1 (\beta,m) &=& \left( \frac{q-1}{1-qy} \right)^{-m+1} \left( \sqrt{y} \right)^{\varepsilon(\beta) -m+1}
 [\beta,m]\\
\noalign{\smallskip}
&=& \left( \frac{q-1}{1-qy} \right)^{-n} \left( \sqrt{y} \right)^{\varepsilon(\alpha) -n+1}
 [\alpha \sigma_n,n+1]\\
\noalign{\smallskip}
&=& \left( \frac{q-1}{1-qy} \right)^{-n} \left( \sqrt{y} \right)^{\varepsilon(\alpha) -n+1} 
\left( \frac{q-1}{1-qy} \right) [\alpha,n]\\
\noalign{\smallskip}
&=& \psi_1(\alpha,n)\,.
\end{array}\]
Suppose that $m=n+1$ and $\beta=\alpha \sigma_n^{-1}$. Observe that the equality $\sigma_n^2 = (q-1) \sigma_n + q$ implies 
\[
\sigma_n^{-1} = q^{-1} \sigma_n - q^{-1} (q-1)\,.
\]
Then
\[\begin{array}{lcl}
\psi_1(\beta,m) &=& \left( \frac{q-1}{1-qy} \right)^{-m+1} \left( \sqrt{y} \right)^{\varepsilon (\beta) -m+1}
[\beta,m]\\
\noalign{\smallskip}
&=& \left( \frac{q-1}{1-qy} \right)^{-n} \left( \sqrt{y} \right)^{\varepsilon (\alpha) -n-1}
[\alpha \sigma_n^{-1},n+1]\\
\noalign{\smallskip}
&=& \left( \frac{q-1}{1-qy} \right)^{-n} \left( \sqrt{y} \right)^{\varepsilon (\alpha) -n-1}
\left( q^{-1} [\alpha \sigma_n, n+1] -q^{-1}(q-1) [\alpha,n+1] \right)\\
\noalign{\smallskip}
&=& \left( \frac{q-1}{1-qy} \right)^{-n} \left( \sqrt{y} \right)^{\varepsilon (\alpha) -n+1} 
\left( \frac{q-1}{1-qy} \right) [\alpha,n]\\
\noalign{\smallskip}
&=& \psi_1(\alpha,n)\,.
\end{array}\]

\bigskip\noindent
By the above, the map $\psi_1$ induces a map $\psi_2: \LL \to \Markov(\sqcup SB)$ defined by  $\psi_2(\hat \beta) = 
\psi_1(\beta)$ for all $\beta \in \sqcup SB$.

\bigskip\noindent
Let $L_+,L_-,L_0$ be three singular links that have the same link diagram except in the neighborhood of a crossing where they 
are like in Figure 1.1.  It is easily deduced from \cite{Birma1} that there exist a singular braid $(\beta,n)$ and an index $1 
\le k \le n-1$ such that $L_+ = \widehat{\beta \sigma_k}$, $L_-=\widehat{\beta \sigma_k^{-1}}$ and $L_0=\widehat{\beta}$. Then 
\[\begin{array}{cl}
&t^{-1} \cdot \psi_2(L_+) -t \cdot \psi_2(L_-)\\
\noalign{\smallskip}
=& \left( \sqrt{yq} \right)^{-1} \left( \frac{q-1}{1-qy} \right)^{-n+1}
\left( \sqrt{y} \right)^{\varepsilon (\beta) -n+2} [\beta \sigma_k,n] - \left( \sqrt{yq} \right) \left(
\frac{q-1}{1-qy} \right)^{-n+1} \left( \sqrt{y} \right)^{\varepsilon (\beta) -n} [\beta \sigma_k^{-1},n]\\
\noalign{\smallskip}
=& \left( \frac{q-1}{1-qy} \right)^{-n+1} \left( \sqrt{y} \right)^{\varepsilon (\beta) -n+1} \left( \frac{1}{\sqrt{q}}
(q-1) [\beta,n] + \frac{1}{\sqrt{q}} q [\beta \sigma_k^{-1},n] - \sqrt{q} [\beta \sigma_k^{-1},n] \right)\\
\noalign{\smallskip}
=& x \left( \frac{q-1}{1-qy} \right)^{-n+1} \left( \sqrt{y} \right)^{\varepsilon (\beta) -n+1} [\beta,n]\\
\noalign{\smallskip}
=& x \cdot \psi_2(L_0)\,.
\end{array}\]
It follows that $\psi_2$ induces a linear map $\psi: \Skein(\LL) \to \Markov (\sqcup SB)$. It is easily checked that this map is 
an algebra homomorphism.

\bigskip\noindent
We turn now to construct the inverse of $\psi$. Let $\phi_1: \sqcup SB \to \Skein (\LL)$ be the map defined by 
\[
\phi_1 (\beta,n) = \left( \frac{q-1}{1-qy} \right)^{n-1} \left( \sqrt{y} \right)^{ n-1 -\varepsilon (\beta)} [\hat \beta]\,.
\]
Let $n \ge 1$, $\alpha, \beta \in SB_n$, and $1\le k\le n-1$.  Then 
\[\begin{array}{cl}
&\phi_1 (\alpha \sigma_k^2 \beta ,n) - (q-1) \cdot \phi_1 (\alpha \sigma_k \beta,n) -q \cdot \phi_1 (\alpha \beta,n)\\
\noalign{\smallskip}
=& \left( \frac{q-1}{1-qy} \right)^{n-1} \left( \sqrt{y} \right)^{ n-1 -\varepsilon (\alpha \beta)} \left( y^{-1}
[ \widehat{ \alpha \sigma_k^2 \beta}] - \left( \sqrt{y} \right)^{-1} 
(q-1) [\widehat{ \alpha \sigma_k \beta}] -q [\widehat{ \alpha \beta}]
\right)\\
\noalign{\smallskip}
=& \left( \frac{q-1}{1-qy} \right)^{n-1} \left( \sqrt{y} \right)^{ n-2 -\varepsilon (\alpha \beta)} \left( \sqrt{q} \right)
\left( t^{-1} [\widehat{ \alpha \sigma_k^2 \beta}] -t [\widehat{\alpha \beta}] -x [\widehat{\alpha \sigma_k \beta}] 
\right)\\
\noalign{\smallskip}
=& 0\,.
\end{array}\]
Let $n\ge 1$ and $\alpha, \beta \in SB_n$. Since $\widehat{\alpha\beta} = \widehat{\beta\alpha}$, we have $\phi_1(\alpha\beta) = \phi_1(\beta\alpha)$. Let $O$ denote the trivial link. One can easily show that
\[
[L \sqcup O] = \left( \frac{t^{-1} -t}{x} \right) [L] = \left( \frac{q-1}{1-qy} \right)^{-1} \left( \sqrt{y} \right)^{-1}
[L] \,,
\]
where $L$ is a link and $L \sqcup O$ is the disjoint union of $L$ and $O$. Now, let $n \ge 1$ and $\alpha \in SB_n$. Observe that $\widehat{(\alpha,n+1)} = \widehat{(\alpha,n)} \sqcup O$, thus 
\[\begin{array}{lcl}
\phi_1 (\alpha, n+1) &=& \left( \frac{q-1}{1-qy} \right)^n \left( \sqrt{y} \right)^{n- \varepsilon (\alpha)} [
\widehat{(\alpha,n+1)}]\\
\noalign{\smallskip}
&=& \left( \frac{q-1}{1-qy} \right)^{n-1} \left( \sqrt{y} \right)^{n-1- \varepsilon (\alpha)} [
\widehat{(\alpha,n)}]\\
\noalign{\smallskip}
&=& \phi_1 (\alpha, n)\,.
\end{array}\]
We also have
\[\begin{array}{lcl}
\phi_1 (\alpha \sigma_n, n+1) &=& \left( \frac{q-1}{1-qy} \right)^n \left( \sqrt{y} \right)^{n- \varepsilon( \alpha
\sigma_n)} [\widehat{( \alpha \sigma_n, n+1)}]\\
\noalign{\smallskip}
&=& z \left( \frac{q-1}{1-qy} \right)^{n-1} \left( \sqrt{y} \right)^{n-1- \varepsilon( \alpha)} 
[\widehat{( \alpha, n)}]\\
\noalign{\smallskip}
&=& z \cdot \phi_1 (\alpha,n)\,.
\end{array}\]
We conclude that the map $\phi_1$ induces a linear map $\phi: \Markov (\sqcup SB) \to \Skein(\LL)$. It is easily checked that $\phi$ is the inverse of $\psi$, thus $\psi$ is an isomorphism.
\qed

\bigskip\noindent
We turn now to state the main result of \cite{ParRab1} from which the calculation of the Markov module of singular braids will 
be deduced.  

\bigskip\noindent
Recall that $S_dB_n$ denotes the set of singular braids with $n$ strands and $d$ singular points, $\HH(S_dB_n)$ denotes the 
subspace of $\HH(SB_n)$ spanned by $S_dB_n$, and that we have the graduation 
\[
\HH (SB_n) = \bigoplus_{d=0}^{+ \infty} \HH (S_dB_n)\,.
\]
Set $\HH_z(S_dB_n) = \C (z,q) \otimes_{\C (q)} \HH (S_dB_n)$.
Let $\Markov(\sqcup S_dB)$ denote the quotient of $\oplus_{n=1}^\infty \HH_z(S_dB_n)$ by the relations
\begin{itemize}
\item
$(ab,n)=(ba,n)$ for all $n\ge 1$ and all $a \in \HH_z(S_kB_n)$ and $b \in \HH_z(S_lB_n)$ such that $k+l=d$;
\item
$(a,n) = (\iota_n(a), n+1)$ for all $n \ge 1$ and all $a \in \HH_z(S_dB_n)$;
\item
$(\iota_n(a) \sigma_n,n+1) = z \cdot (a,n)$ for all $n \ge 1$ and all $a \in \HH_z(S_dB_n)$.
\end{itemize}
It is clear that
\[
\Markov (\sqcup SB) = \bigoplus_{d=0}^{+ \infty} \Markov (\sqcup S_dB)\,.
\]

\bigskip\noindent
Let $\C(q,z)[S_dB_n]$ be the vector space over $\C (q,z)$ freely spanned by $S_dB_n$. For $d \ge 1$, we define the linear maps 
$f_{n,0},f_{n,1} : \C(q,z)[S_dB_n] \to \C(q,z)[S_{d-1}B_n]$ as follows. Let $\beta \in S_dB_n$. Then write $\beta$ in the form 
$\beta= \alpha_0 \tau_{i_1} \alpha_1 \cdots \tau_{i_d} \alpha_d$ with $\alpha_i \in B_n$ for $0 \le i\le d$, and set
\begin{align*}
f_{n,0} (\beta) &= \sum_{k=0}^d \alpha_0 \tau_{i_1} \alpha_1 \cdots \tau_{i_{k-1}} \alpha_{k-1} \alpha_k \tau_{i_{k+1}}
\alpha_{k+1} \cdots \tau_{i_d} \alpha_d\\
f_{n,1} (\beta) &= \sum_{k=0}^d \alpha_0 \tau_{i_1} \alpha_1 \cdots \tau_{i_{k-1}} \alpha_{k-1} \sigma_{i_k}
\alpha_k \tau_{i_{k+1}} \alpha_{k+1} \cdots \tau_{i_d} \alpha_d
\end{align*}
It follows from Theorem 2.1 that this definition does not depend on the choice of the expression of $\beta$.

\bigskip\noindent
It is easily checked that the collection of linear maps $\{f_{n,0}\}_{n\ge 1}$ induces a linear map $g_0:\Markov(\sqcup S_dB) \to 
\Markov(\sqcup S_{d-1}B)$. Similarly, the collection of maps  $\{f_{n,1}\}_{n\ge 1}$ induces a linear map $g_1:\Markov(\sqcup S_dB) 
\to \Markov(\sqcup S_{d-1}B)$.

\bigskip\noindent
Let $\Markov(\sqcup S_dB)^\ast$ be the dual space of $\Markov(\sqcup S_dB)$, that is, the space of linear forms on 
$\Markov(\sqcup S_dB)$. For $d\ge 1$, we denote by $\Phi_{d,0}:  \Markov(\sqcup S_{d-1}B)^\ast \to  \Markov(\sqcup S_dB)^\ast$ 
the linear map induced by $g_0$, and by $\Phi_{d,1}:  \Markov(\sqcup S_{d-1}B)^\ast \to  \Markov(\sqcup S_dB)^\ast$ the linear 
map induced by $g_1$. Note that $\Phi_{d+1,1} \circ \Phi_{d,0} = \Phi_{d+1,0} \circ \Phi_{d,1}$ for all $d \ge 1$.

\bigskip\noindent
For  $d \ge 0$, we define the elements $T_{d,0}, T_{d,1}, \dots, T_{d,d} \in \Markov(\sqcup S_dB)^\ast$ by induction on $d$ as 
follows. It is proved in \cite{Jones1} that the space $\Markov(\sqcup S_0B)^\ast$ is of dimension $1$. Then we denote by $T_{0,0}$ the generator of 
$\Markov(\sqcup S_0B)^\ast$ whose value on the trivial braid is $1$. Suppose $d \ge 1$. Then we set
\[\begin{array}{lcll}
T_{d,0} &=& \Phi_{d,0} (T_{d-1,0})\\
T_{d,k} &=& \Phi_{d,0} (T_{d-1,k}) = \Phi_{d,1} (T_{d-1,k-1}) &\quad \text{if } 1 \le k\le d-1\\
T_{d,d} &=& \Phi_{d,1} (T_{d-1,d-1})
\end{array}\]

\bigskip\noindent
{\bf Theorem 2.7} (Paris, Rabenda \cite{ParRab1}).
{\it Let $d \ge 0$. Then $\Markov(\sqcup S_dB)^\ast$ is of dimension $d+1$, and $\{T_{d,0}, T_{d,1}, \dots, T_{d,d}\}$ is a 
basis for $\Markov(\sqcup S_dB)^\ast$.}

\bigskip\noindent
We can now calculate the Markov module of singular braids:

\bigskip\noindent
{\bf Theorem 2.8.}
{\it The algebra $\Markov(\sqcup SB)$ is a polynomial algebra $\C(q,z)[X,Y]$ in two variables $X$ and $Y$, where $X$ and $Y$ 
are the classes of $\tau_1$ and $\tau_1\sigma_1$, respectively.}

\bigskip\noindent
{\bf Proof.}
Let  $d \ge 0$ and $0 \le k\le d$. Observe that $X^kY^{d-k}$ is the class of $\tau_1 \cdots \tau_{2k-1} (\tau_{2k+1} 
\sigma_{2k+1})
\break
\cdots (\tau_{2d-1} \sigma_{2d-1})$. In particular, we have $X^kY^{d-k} \in \Markov(\sqcup S_dB)$. So, in order 
to prove Theorem 2.8, it suffices to show that $\{ X^d, X^{d-1}Y, \dots, XY^{d-1},Y^d\}$ is a basis for $\Markov(\sqcup S_dB)$. 
Since we already know by Theorem 2.7 that $\Markov(\sqcup S_dB)$ is of dimension $d+1$, it actually suffices to show that $\{ 
X^d, X^{d-1}Y, \dots, XY^{d-1},Y^d\}$ is linearly independent. We prove this by induction on $d$. The case $d=0$ being trivial, 
we assume $d\ge 1$ plus the inductive hypothesis.

\bigskip\noindent
A direct calculation shows that, for $0 \le k\le d$, we have
\[\begin{array}{lcl}
g_0( X^k Y^{d-k}) &=& k \cdot X^{k-1} Y^{d-k} +z(d-k) \cdot X^k Y^{d-k-1}\\
\noalign{\smallskip}
g_1( X^k Y^{d-k}) &=& kz \cdot X^{k-1} Y^{d-k} + (d-k) ((q-1)z+q) \cdot X^k Y^{d-k-1}
\end{array}\]
Let $a_0, a_1, \dots, a_d \in \C(q,z)$ such that
\begin{equation}
\sum_{k=0}^d a_k X^k Y^{d-k} =0\,.
\label{eq22}
\end{equation}
Applying $g_0$ and $g_1$ to (2.2) we obtain
\begin{gather*}
\sum_{k=0}^{d-1} \big( (k+1) a_{k+1} + (d-k)za_k \big) X^k Y^{d-k-1} =0\\
\sum_{k=0}^{d-1} \big( (k+1)za_{k+1} + (d-k)((q-1)z+q)a_k \big) X^k Y^{d-k-1} = 0
\end{gather*}
By induction, it follows that
\[
\left\{\begin{array}{rcrcl}
(k+1)\, a_{k+1} &+& (d-k)z\, a_k &=& 0\\
(k+1)z\, a_{k+1} &+& (d-k)((q-1)z+q)\, a_k &=& 0
\end{array}\right.
\]
for all $0\le k\le d-1$. The determinant of this system of linear equations in the variables $a_{k+1},a_k$ is equal to 
$-(k+1)(d-k)(z^2-(q-1)z-q) \neq 0$, thus $a_{k+1}=a_k=0$.
\qed

\bigskip\noindent
{\bf Corollary 2.9.}
{\it Let $\LL$ be the set of singular links in the sphere $\S^3$. Then the algebra $\Skein(\LL)$ is a polynomial algebra 
$\C(t,x)[\hat X, \hat Y]$ in two variables $\hat X$ and $\hat Y$, where $\hat X$ and $\hat Y$ are the classes of the singular 
links $L_X$ and $L_Y$ represented in Figure 1.4.}

\bigskip\noindent
{\bf Question 2.10.}
The proof that the set $\BB = \{X^aY^b; a,b \in \N\}$ is linearly independent in $\Markov (\sqcup SB)$ is entirely given in the 
above proof  of Theorem 2.8, and does not need Theorem~2.7 at all. However, the proof that $\BB$ spans $\Markov (\sqcup SB)$ 
uses the fact that the dimension of  $\Markov (\sqcup S_dB)$ is (less or) equal to $d+1$ for all $d \ge 0$, and the proof of 
this latest assertion needs long and tedious calculations. It would be interesting to find a (simplest) proof of the equivalent 
fact that $\hat \BB = \{ \hat X^a \hat Y^b; a,b \in \N\}$ spans $\Skein(\LL)$, which would directly use the Skein relations.



\bigskip\bigskip\noindent
{\bf Luis Paris},

\smallskip\noindent 
Institut de Math\'ematiques de Bourgogne, UMR 5584 du CNRS, Universit\'e de Bourgogne, B.P. 
47870, 21078 Dijon cedex, France

\smallskip\noindent
E-mail: {\tt lparis@u-bourgogne.fr}


\end{document}